\documentclass[12pt,a4paper]{amsart}
\usepackage{graphicx} 
\usepackage[a4paper, left=2.5cm, right=2.5cm, top=2.8cm, bottom=2.8cm]{geometry}
\usepackage{natbib}
\usepackage{hyperref}
\usepackage{xcolor}
\usepackage{amsaddr}

\newcommand\Ex{\text{\sf E}}
\newcommand\Prob{\text{\sf P}}
\newcommand\Var{\text{\sf V}}
\newcommand\Cov{\text{\sf Cov}}

\newtheorem{theorem}{Theorem}
\newtheorem{proposition}{Proposition}
\newtheorem{lemma}{Lemma}
\newtheorem{corollary}{Corollary}
\newtheorem{definition}{Definition}
\theoremstyle{remark}

\newtheorem{remark}{Remark}
\begin{document}
\title[Asymmetric switch process]{Characteristics of asymmetric switch processes with independent switching times}

\author[H. Bengtsson and K. Podg\'orski]{Henrik Bengtsson and Krzysztof Podg\'orski} \address{Department of Statistics,
Lund University, \\}
\email{ Henrik.Bengtsson@stat.lu.se, \href{https://orcid.org/0000-0002-9280-4430}{ \includegraphics[height=2.2mm]{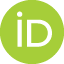} 0000-0002-9280-4430}}  
\email{ krzysztof.podgorski@stat.lu.se, \href{https://orcid.org/0000-0003-0043-1532}{ \includegraphics[height=2.2mm]{ORCIDiD64.png} 0000-0003-0043-1532}} 

\date{January 27, 2025}
\begin{abstract}
The asymmetric switch process is a binary stochastic process that alternates between the values one and minus one, where the distributions of the time in these states may differ. Two versions of the process are considered: a non-stationary version that starts with a switch at time zero and a stationary one constructed from the non-stationary one. Characteristics of these two processes, such as the expected values and covariance, are investigated. The main results show an equivalence between the monotonicity of the expected value functions and the distribution of the intervals having a stochastic representation in the form of a sum of random variables, where the number of terms follows a geometric distribution. This representation has a natural interpretation as a model in which switching attempts may fail at random. From these results, conditions are derived when these characteristics lead to valid interval distributions, which is vital in applications. 
\end{abstract}

\maketitle


\section{Introduction} 
\noindent 
The study of binary stochastic processes has a long and rich history. Throughout this history, various particular forms have received considerable attention, from elementary Markov chains to renewal processes. The attention garnered by these processes is understandable due to their many applications in fields such as queuing theory \cite{Erlang}, signal processing \citep{SBRS}, statistical physics \citep{IIA}, and approximating the excursion distribution for stochastic processes \citep{McF56} and \cite{BK}.

To study the excursion of stochastic processes \cite{McF56} introduced the notion of the clipped process. The clipped process is obtained from a sufficiently smooth process $W(t)$ by computing $Z(t)={\rm sign} (W(t)-u)$. The process $Z(t)$ is a binary process taking the values one and minus one, and the time spent in the two states represents the excursion times of $W(t)$, i.e., the times $W(t)$ is above or below the value $u$. Connecting characteristics, such as the covariance function of $Z(t)$ to the excursion distributions of $W(t)$ is still a notoriously tricky problem, despite having received considerable attention during the 19th century by, among others, \cite{Rice} and \cite{McF56, McF58}.

The difficulty stems from the complex dependency structures of the intervals of the clipped process $Z(t)$. One way to address this difficulty is through approximation methods. One such method is the independent interval approximation (IIA), and for an overview of this method, see \cite{Sire2007, Sire2008} and \cite{IIA} for extensive overviews of this method in the context of applications in physics. This method uses a switch process to approximate $Z(t)$ and, through this, approximate the distributions of the $u-$level excursions of $W(t)$.

The switch process is a binary process similar to the clipped process. However, the times the process takes the value one and minus one are independent. This independence significantly reduces the complexity of linking the process's characteristics to the distribution of the intervals. Thus, characteristics of $Z(t)$ can be imposed on a switch process, and then the distributions of the interval lengths are deduced from these. This deduced distribution then serves as an approximation for the excursion distribution of $W(t)$.

This application motivates the study of the relationship between the switch process's characteristics and the distribution of the intervals. A partial characterization was obtained in \cite{Bengtsson2} for when the interval lengths have the same distribution for both states. In this paper, we allow these interval lengths to have different distributions. In this sense, the switch process is called asymmetric, in contrast to the symmetric one studied in \cite{Bengtsson2}. This extension has implications for the previously mentioned IIA framework when approximating non-zero-level excursions.

The outline of the paper is as follows. The two versions of the switch process are introduced in Section \ref{sec:SP}. The main results are found in Section \ref{sec:main}. The paper concludes with some illustrating examples in Section \ref{sec:ex}. Some supporting results and proofs have been collected in the Appendix to keep the paper somewhat self-contained.

\section{Switch process} \label{sec:SP}
\noindent
Two versions of the switch process are constructed from two independent switching time distributions, which are the distributions of the length of the intervals that the process spends in each state. Let $T_+$ and $T_-$ denote the random variables associated with these distributions. Where $T_+$ is the length of the intervals for which the process takes the value one, and $T_-$ is the length of the intervals for which the process takes the value minus one. While the details on how they are used for constructing the switch processes will be discussed later, the two versions are illustrated in Figure \ref{SwitchLine} and \ref{DelSw}.

For the construction of the processes, we need some basic assumptions on the distributions of $T_+$ and $T_-$. Unless otherwise stated, it will be assumed that they are non-negative, have no atom at zero, are absolutely continuous with support on the entire non-negative real line, and have finite expectation. Additionally, we assume that there exists a $m\in \mathbb{N}: \sup_{t>0} (f_+ \ast f_-)^{\ast m}(t)<\infty$, where $f_+$ and $f_-$ are the densities associated with $T_+$ and $T_-$ and $\ast$ denotes convolution of functions. This condition is necessary for the existence of the derivatives of some characteristics, such as the expected value functions.

To study these functions, we use the Laplace transform, which is the standard tool in renewal theory. We denote the Laplace transform by $\mathcal{L}(\cdot)$ and in particular the Laplace we denote the Laplace transform of a probability  density by $\Psi(s)$. Some elementary results and the definition can be found in the appendix. With the notation of the Laplace transform in place, we are now ready to define the non-stationary switch process.


\subsection*{The non-stationary switch process}
\begin{figure}
\hspace{-15mm}
\includegraphics[width=0.75\textwidth]{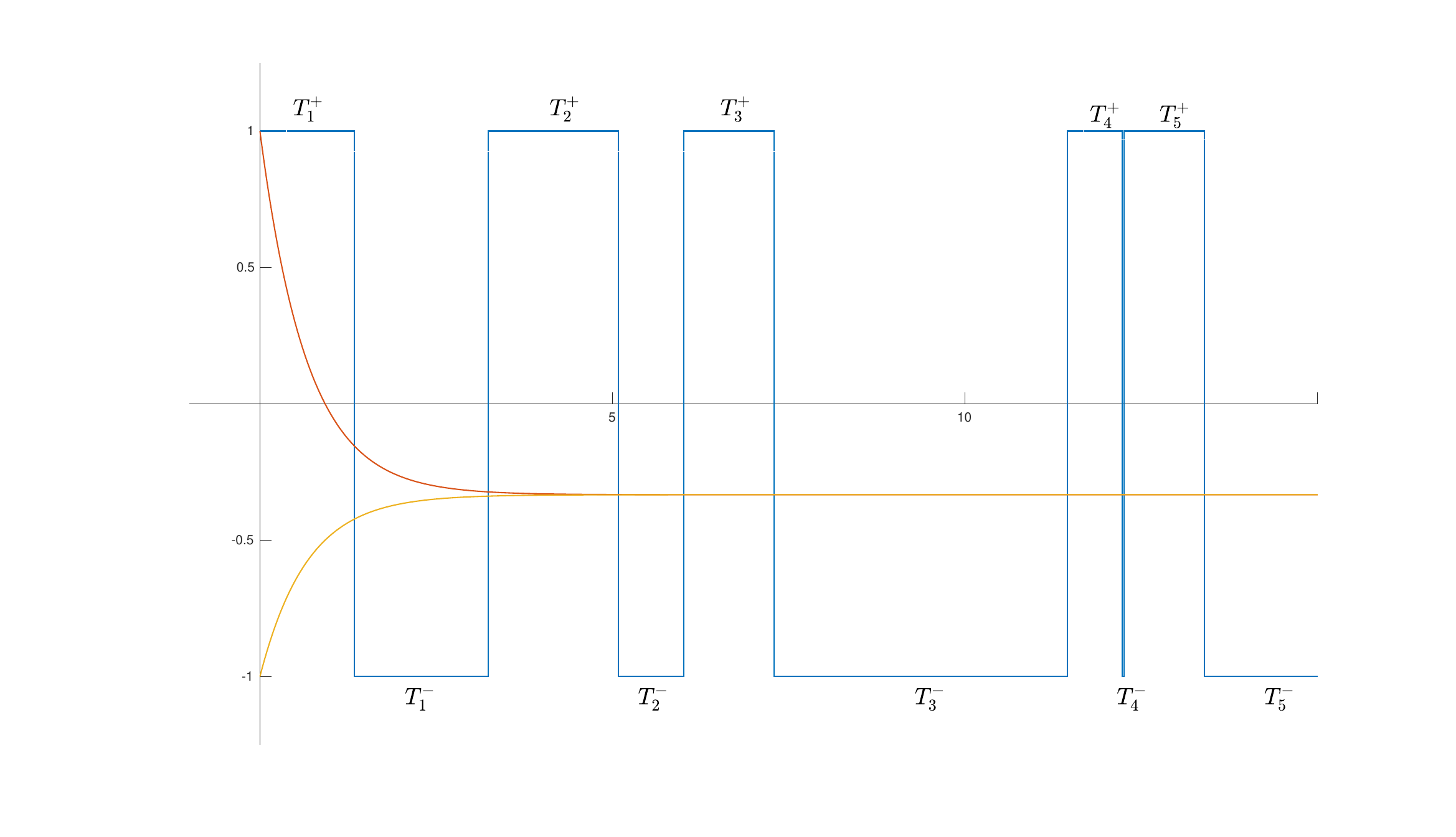}\vspace{-5mm}
\caption{
A realization of a non-stationary switch process (blue), $D(t)$, with exponential switching times and $\delta=1$. Together with $E_+(t)$ (red) and $E_-(t)$ (orange).
}
\label{SwitchLine}
\end{figure}
The construction of the non-stationary switch process starts with the binary random variable $\delta$, independent of the switching time distributions, with $\Prob{(\delta=1)}=p$ and $\Prob{(\delta=-1)}=1-p$. The two values of $\delta$ are denoted by $+$ and $-$ when $\delta$ is used as a subscript. This random variable $\delta$ determines the starting value of the process at the origin. We then interlace random intervals of lengths distributed as  $T_+$ and $T_-$, respectively, starting from the origin. Since we now consider sequences of switching time distributions, we denote these by $T_i^+$, $i \in \mathbb N$ and $T_i^-$, $i \in \mathbb N$. If $\delta=1$, the non-stationary switch process takes the value one over the interval $[0, T_1^+)$, and then switches to minus one on the interval $[T_1^+, T_1^++T_1^-)$, then it switches back to one and so forth. If $\delta=-1$, the process takes the value minus one on the interval $[0, T_1^-)$ and then it switches to one on $[T_1^-, T_1^-+T_1^+)$ and so on when time passes. The role $\delta$, $T_+$, and $T_-$ play in the construction of the non-stationary switch process is perhaps best illustrated in Figure \ref{SwitchLine}, which shows a realization of the process.

We have so far not placed any assumptions on the dependency structure of $T_i^+$, $i \in \mathbb N$, $T_i^-$, $i \in \mathbb N$. However, we only consider the simple cases where $T_i^+$ and $T_i^-$ constitute two mutually independent sequences of independent identically distributed (iid) random variables. A binary stochastic process taking the values $\{-1,1\}$ and is constructed under the mentioned assumption by interlacing intervals, as described previously, is called a non-stationary switch process and is denoted by $D(t)$, $t\geq0$.  

By having a switch at zero, the process will not be stationary in the traditional sense; it will, however, be cycle stationary. This notion of stationarity is defined by the cycles $(T_i^+,T_i^-)$, $i\in \mathbb{N}$, having the same distribution, i.e ${(T_i^+,T_i^-)\stackrel{d}{=}(T_j^+,T_j^-)}$ for all $i,j\in \mathbb{N}$, where $\stackrel{d}{=}$ denotes equality in distribution. It should be noted that cycle stationery will be preserved by letting $T_i^+$ and $T_i^-$ be dependent. However, this case will not be investigated further in this paper.

The lack of stationary in the traditional sense makes the covariance function depend on two arguments, which is harder to utilize. This leads us to investigate the expected value function conditioned on the starting value $\delta$. The relation between the expected value functions and the distributions of $T_+$ and $T_-$ in the Laplace domain is given next.
\begin{proposition}
\label{OneDim}
The Laplace transform of $P_\delta(t)=\Prob{\left(D(t)=1 \big | \delta\right)}$, $t>0$, is given by
\begin{align*}
\mathcal L (P_\delta)(s)&=\frac{1-\Psi_{+}(s)}{s(1-\Psi_+(s)\Psi_-(s))} \cdot
\begin{cases}
1&;\delta=1,\\
\Psi_-(s)&; \delta=-1,
\end{cases}
\end{align*}
where $\Psi_+$ and $\Psi_-$ are Laplace transformation of probability distributions corresponding to $T_+$ and $T_-$, respectively. 
Moreover, for the expected value function $E_\delta(t)=\Ex{( D(t)\big |\delta)}$, $t>0$ we have 
\begin{align*}
\mathcal L (E_\delta(s))&= \frac{\Psi_-(s)-\Psi_+(s)+\delta (1-\Psi_-(s))(1-\Psi_+(s))}{s(1-\Psi_+(s)\Psi_-(s))}.
\end{align*}  
\end{proposition}
In addition to this proposition, we have the following corollary on the derivatives of $E_+$ and $E_-$, which will be important when investigating switch processes with monotone expected value functions. 

\begin{corollary}\label{cor:Inv}
    Let $E_+$ and $E_-$ be the expected value functions of a non-stationary switch process with switching time distributions such that there exists a $m\in \mathbb{N}$: \\$\sup_{t>0} (f_+ \ast f_-)^{\ast m}(t)<\infty$. Then we have the following relations
    \begin{align*}
        \mathcal L (E_+')(s)=-2\Psi_+(s)\frac{1-\Psi_-(s)}{1-\Psi_+(s)\Psi_-(s)},&\,\,\, \mathcal L (E_-')(s)=2\Psi_-(s)\frac{1-\Psi_+(s)}{1-\Psi_+(s)\Psi_-(s)}.
    \end{align*}
\end{corollary}
This proposition and corollary focuses only on the expected value functions, its derivatives, and $P_\delta(t)$. More general distributional properties are difficult to derive. For example, obtaining an explicit form for the covariance function is difficult. A more straightforward function to derive is the variance of $D(t)$. Since $\Ex(D(t)^2)=1$, we have
\begin{align*}
    \Var (D)(t)= 1-\left(\Ex D(t) \right)^2.
\end{align*}
From the above equation, it is clear that the variance diminishes closer to zero since the probability that a switch has occurred is fairly small.

The limiting behavior of the expected value functions is important for describing the process's behavior and characterizing the expected value functions. While these limits follow standard results on alternating renewal processes, we present a rigorous treatment in Lemma \ref{lem:Ps}, which can be found in the appendix. From this lemma and the observation that $E_\pm(t)=2P_\pm(t)-1$ we have the following limits 
\begin{align} \label{eq:lim}
    \begin{split}
        &\lim_{t\rightarrow 0^+} E_+(t)=1, \ \ \ 
        \lim_{t\rightarrow 0^+} E_-(t)=-1,
        \\
        &\lim_{t\rightarrow \infty} E_+(t)=\lim_{t\rightarrow \infty}E_-(t)= \frac{\mu_+-\mu_-}{\mu_++\mu_-}.
    \end{split}
\end{align}
From these limits, it is clear that the initial effect of the switch placed at zero will diminish when $t$ goes to infinity.

Until now, we have only considered the non-stationary switch process on the non-negative part of the real line. However, the process can be extended to the entire real line, similar to how it was constructed for the non-negative part. Suppose $\delta=1$ then the extended process would take the value minus one on the interval $[-T_{-1}^{-},0)$ and then switch to one on $[-(T_{-1}^++T_{-1}^-), -T_{-1}^-)$ and so forth. If $\delta=-1$ the process becomes one on  $[-T_{-1}^{+},0)$ and minus one on $[-(T_{-1}^++T_{-1}^-), -T_{-1}^+)$. Hence, for a given sequence of pairs $\mathcal T= (T_i^-, T_i^+)_{i \in \pm \mathbb N}$, and initial choice of the sign $\delta$, can be used to extend the non-stationary switch process on the entire real line. 

The attachment of a switch at zero results in the process not being stationary in the traditional sense. In the subsequent subsection, a standard technique from renewal theory is used to construct a stationary version by delaying the process backward and forward around zero.

\subsection*{The stationary switch process}
\begin{figure}[!t]
\vspace{-15mm}
\includegraphics[width=0.75\textwidth]{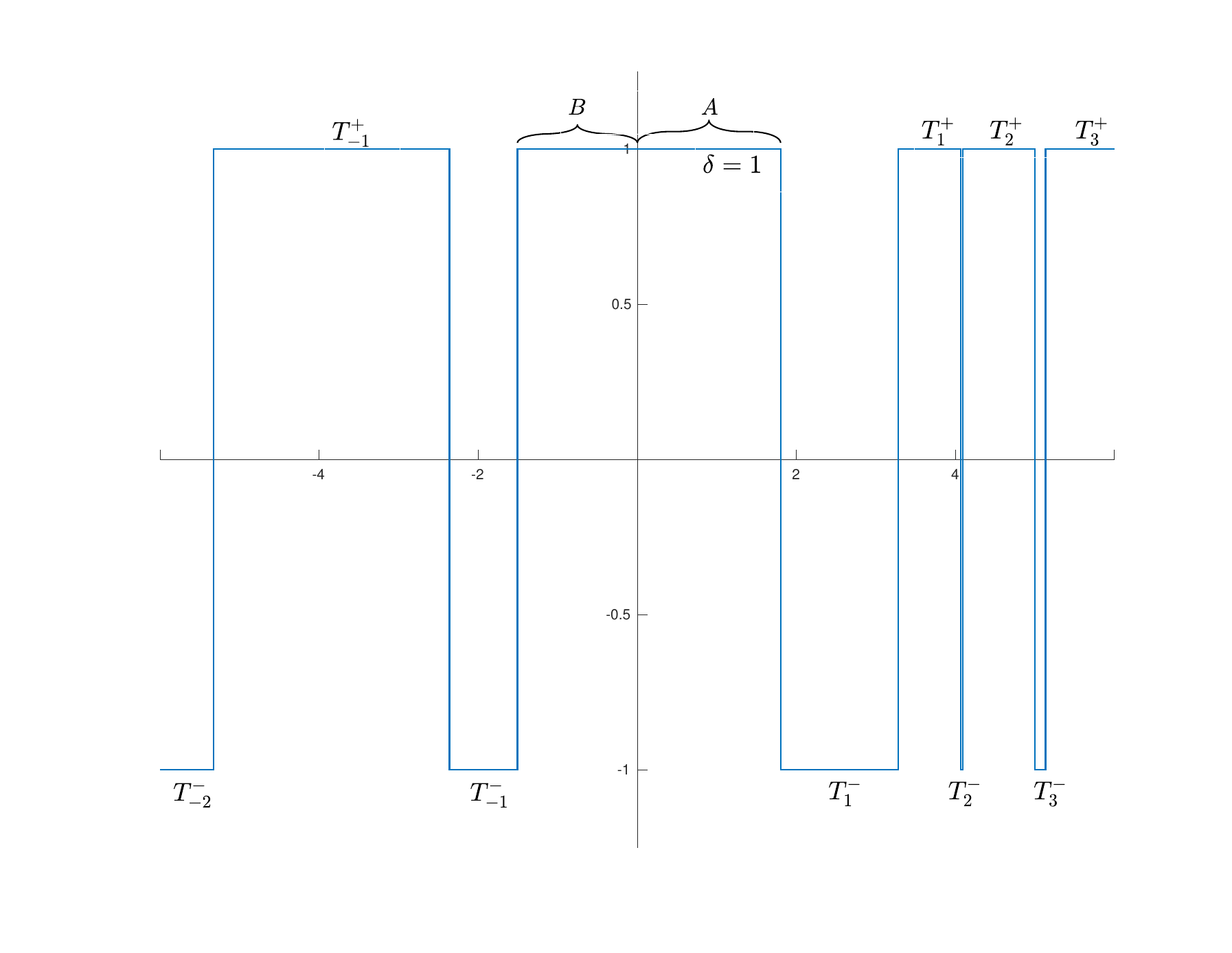}
\vspace{-10mm}
\caption{ 
A realization of a stationary switch process, $\tilde{D}$ with exponential switching times and $\delta=1$. $A$ and $B$ are the forward and backward delays. }
\label{DelSw}
\end{figure}
The stationary switch process is constructed from the non-stationary switch process. Since this process is cycle stationary, it follows from \cite{Herman} that such a construction is possible. This construction is done by delaying the process forward and backward so there is no longer a switch at the origin. Hence, the interval containing the origin will follow a different distribution than $T_+$ or $T_-$. In essence, we reintroduce the well-known inspection paradox for which the mathematical details have been known since \cite{Palm}.

For the construction, we follow a similar argument as in \cite{Joint22}. Let $[-B, A]$ be the interval containing the origin, and $\delta$ be the value of the process on this interval. From the ends, $-B$ and $A$ non-stationary switch processes are attached with the initial values $-\delta$. The process with these delays is called the delayed switch process and is denoted by $\tilde{D}(t)$, $t \in \mathbb{R}$. A realization of $\tilde{D}(t)$ can be seen in Figure \ref{DelSw}, and the problem now becomes finding the distribution of $(A, B, \delta)$ such that the delayed process becomes stationary. The solution is presented in the next proposition and follows from a standard application of the key renewal theorem, which can be found in the appendix. 
\begin{proposition}
\label{stat}
Let $\mu_+=E T_+$, $\mu_-=E T_-$ and $f_+$, $f_-$ be the densities of $T_+$ and $T_-$, respectively. 
If the distribution of $(A,B,\delta)$ is given through
\begin{align*}
\Prob{(\delta = 1)}=\frac{\mu_+}{\mu_-+\mu_+},&~
\Prob{(\delta = -1)}=\frac{\mu_-}{\mu_-+\mu_+}\\
f_{A,B|\delta}(a,b|1)=\frac{f_+(a+b) }{\mu_+}, &~ f_{A,B|\delta}(a,b|-1)=\frac{f_-(a+b) }{\mu_-},
\end{align*}
where $f_{A, B|\delta}$ stands for the conditional density, then the delayed switch process $\tilde D(t)$ is stationary in time, with the expected value $(\mu_+-\mu_-)/(\mu_++\mu_-)$. 
Additionally its distribution is uniquely characterized by $\tilde P_\delta(t)=\Prob(\tilde D(t)=1|\tilde D(0)=\delta)$ that has the Laplace transform of the following form
\begin{align*}
\mathcal L \tilde P_\delta(s)
&=
\frac{1}{s}
\begin{cases} \displaystyle
1-\frac 1 {\mu_+ s} 
\frac{\left(1-\Psi_+(s)\right)\left(1-\Psi_-(s)\right)}{1-\Psi_+(s)\Psi_-(s)} 
&;\delta=1,\vspace{2mm}\\
\displaystyle
\frac 1 {\mu_- s} 
\frac{\left(1-\Psi_+(s)\right)\left(1-\Psi_-(s)\right)}{1-\Psi_+(s)\Psi_-(s)}&; \delta=-1,
\end{cases}
\end{align*}
where $\Psi_+$ and $\Psi_-$ are Laplace transformation of probability distributions corresponding to $T_+$ and $T_-$, respectively. 
The autocovariance function for this process, $R(t)=\Cov{(\tilde D(u), \tilde D(u+t))}$, is on the form
\begin{align*}
    R(t)=\frac{2}{\mu_++\mu_-}
    \left(
    \tilde P_{+}(t)\mu_+-\tilde P_{-}(t)\mu_-+
    \mu_+ 
    \frac{\mu_--\mu_+}{\mu_++\mu_-}
    \right).
\end{align*}
Lastly, the Laplace transform of $R(t)$, $t>0$, is given by
\begin{align*}
    \mathcal L (R)(s)=\frac{4}{s\left(\mu_++\mu_-\right)}\left(\frac{\mu_+\mu_-}{\mu_++\mu_-}- \frac1s \frac{(1-\Psi_+(s))(1-\Psi_-(s))}{1-\Psi_-(s)\Psi_+(s)} \right) 
\end{align*}
\end{proposition}
A delayed switch process with $[-B, A]$ and $\delta$ as described in Proposition \ref{stat} will, from this point onwards, be called a stationary switch process and denoted by $\tilde{D}(t)$. While we assume that the densities of $T_+$ and $T_-$ exist, this is not a technical necessity, as seen in the following remark. 

\begin{remark}
\label{AB}
We should note from the proof that it is not required for the distribution of the switching times to have a density. 
It is sufficient that the distribution of $T_1^++T_2^-$ is not sitting on a lattice for the key renewal theorem to hold. In fact, the stationary distribution of $(A,B,\delta)$ is given by
\begin{align*}
\Prob{(\delta=1, A\le a, B>b)}=\frac{\int_0^a \bar F_+(u+b)~du}{\mu_++\mu_-}
\end{align*}
We note two simple consequences of this fact. 
First, the distributions of the delays $A$ and $B$ are identical and always continuous with respect to the Lebesgue measure with their densities conditionally on $\delta$ given by 
\begin{align*}
f_{A|\delta}(x|\delta)=f_{B|\delta}(x|\delta) =\frac{\bar F_\delta(x)}{\mu_\delta},
\end{align*}
which in the Laplace domain becomes 
\begin{align*}
\Psi_{A|\delta}(s)=\Psi_{B|\delta}(s)=\frac{1-s\mathcal L F_\delta(s)}{s\mu_\delta}=\frac{1-\Psi_\delta(s)}{s\mu_\delta}.
\end{align*}
using Equation (\ref{cdf}).
Second, we note that the distribution of the interval containing the origin is given by the density
\begin{align*}
 f_{A+B|\delta}(x|\delta)=\frac{x f_\delta(x)}{\mu_\delta}.   
\end{align*}
\end{remark}

We have now constructed a stationary switch process $\tilde{D}(t)$ from the non-stationary switch process $D(t)$. In the next section, we connect the characterizing function to the properties of the switching time distributions.

\section{Characteristics of the switch process} \label{sec:main}
\noindent
We have seen how the distributions of $T_+$ and $T_-$ give rise to the following functional characteristics: the expected value functions $E_+$, $E_-$, and $P_+$, $P_-$ of the non-stationary version, and the autocovariance $R$, $\tilde P_+$, and $\tilde P_-$ for the stationary counterpart. Now, we reverse the question: What functional characteristics lead to valid distributions, and what are their properties? To address these questions, we start with a proposition showing that the distributions of $T_+$ and $T_-$ are identifiable from the expected value functions of the non-stationary switch process.
\begin{proposition} \label{prop:invert}
    The distributions $T_+$ and $T_-$ are uniquely identifiable from $E_+$ and $E_-$ in the Laplace domain through the equations 
\begin{align*} 
        \Psi_+(s)&=\frac{s\mathcal{L}(E_+)(s)-1}{s\mathcal{L}(E_-)(s)-1},  
        \\
        \Psi_-(s)&=\frac{s\mathcal{L}(E_-)(s)+1}{s\mathcal{L}(E_+)(s)+1}.
\end{align*}
\end{proposition}

Proposition \ref{prop:invert} shows that the distribution of $T_+$ and $T_-$ can be uniquely determined from the expected value functions. The subsequent proposition provides the basis to examine if this is also possible from the covariance function.
\begin{proposition}
\label{prop:rel}
Using the definitions of Proposition~\ref{OneDim} and \ref{stat}, we have the following relations between characteristics of the stationary and non-stationary switch process:
\begin{align*}
    &\mathcal L(\tilde P'_+)=\frac{\mathcal L(P_-) -\mathcal L(P_+)}{\mu_+},\,\, 
    \mathcal L(\tilde P'_-)= \frac{\mathcal L(P_+) -\mathcal L(P_-)}{\mu_-}. \\
    &\mathcal L(R')=4\,\frac{\mathcal L(P_-)-\mathcal L(P_+)}{\mu_++\mu_-}=2\,\frac{\mathcal L(E_-)-\mathcal L(E_+)}{\mu_++\mu_-}.
\end{align*}
\end{proposition}

The following remark presents the implications of this proposition for the recoverability of $T_+$ and $T_-$ from the covariance function. 
\begin{remark}
    The probabilities $\tilde{P_+}$ and $\tilde{P_-}$ can not be uniquely obtained from $P_+$ and $P_-$. Subsequently any relationship between $P_+$, $P_-$ and $T_+$, $T_-$ can not be used to derive $\tilde{P_+}$ and $\tilde{P_-}$. This lack of identifiability is shown by considering the following system of equations, which follows directly from the top equations of Proposition \ref{prop:rel}
    \begin{align*}
   \begin{bmatrix}
    \mathcal L(\tilde P'_+)(s) \\
    \mathcal L(\tilde P'_-)(s) 
    \end{bmatrix}
    =
    \begin{bmatrix}
    -\frac{1}{\mu_+} & \frac{1}{\mu_+}\\
    \frac{1}{\mu_-} & -\frac{1}{\mu_-}
    \end{bmatrix}
    \begin{bmatrix}
    \mathcal L (P'_+)(s) \\
    \mathcal L (P'_-)(s)
    \end{bmatrix}.
\end{align*}
Since the above $2 \times 2$ matrix is singular, there is no unique solution such that $\mathcal L (P'_+)$ and $\mathcal L( P'_-)$ can expressed by $\mathcal L (\tilde{P'_+})$ and $\mathcal L (\tilde{P'_-})$. 
\end{remark}
A consequence of the previous remark is that if one wants to recover $T_+$ and $T_-$ from characteristics of the stationary switch process, it cannot be done purely from $R$, $\tilde{P_+}$ and $\tilde{P_+}$. Thus, additional characteristics need to be considered. One might consider the expected value of the number of switches from one to minus one and vice versa. However, these will not provide additional information to identify the distributions of $T_+$ and $T_-$ since they are linear in time due to the stationarity of $\tilde{D}$. Therefore, it is a nontrivial task to retrieve the switching time distributions from the characteristics of the stationary switch process.

In contrast to the previously mentioned difficulties, Proposition \ref{prop:invert} seems to provide a straightforward way to recover the switching time distributions. However, this is not the case. For a proposed pair of expected value functions, it has to be verified that the right-hand side of the two equations of Proposition \ref{prop:invert} are completely monotone functions. Recall that a function $\phi$ on $[0,\infty)$ is completely monotone if it has finite limits, is infinitely differentiable, and satisfies the condition for all $n\in \mathbb{N}$
\begin{align*}
    (-1)^n \phi^{(n)}(s)\ge 0, \,\, s>0.
\end{align*}
Hence, to show that the two equations of Proposition \ref{prop:invert} are completely monotone is arduous at best and almost impossible at worst. This difficulty is a consequence of Bernstein’s theorem (see Theorem 1, page 439 in \cite{FellerV2}), which states that a function is the Laplace transform of a probability distribution if and only if it is completely monotone and one at zero. Therefore, recovering the switching time distributions from Proposition \ref{prop:invert} remains a problem. While the following proposition does not resolve the issue, it characterizes the properties of the expected value functions. 

\begin{proposition}
\label{prop:char}
The following conditions are necessary and sufficient for differentiable functions $E_+$ and $E_-$ on $[0,\infty)$ with values in $[-1,1]$ to define a valid switch process $D(t)$ for which they correspond to $\Ex\left(D(t)|D(0)=1\right)$ and $\Ex\left(D(t)|D(0)=-1\right)$, respectively,
\begin{itemize}
 \item[\it i)] $\displaystyle \lim_{t\rightarrow 0^+} E_+(t)=1$, $\displaystyle \lim_{t\rightarrow 0^+} E_-(t)=-1$, $\displaystyle \lim_{t\rightarrow \infty} E_+(t)=\lim_{t\rightarrow \infty} E_-(t)=\gamma$,  for some $\gamma\in (-1,1)$, \vspace{2mm}
 \item[\it ii)] the functions $\displaystyle\frac{\mathcal L (E_+')(s)}{\mathcal L (E_-')(s)-2}$ and $\displaystyle\frac{\mathcal L (E_-')(s)}{\mathcal L (E_+')(s)+2}$ are completely monotone.
\end{itemize}
\end{proposition}

For the class of non-stationary switch processes with monotone expected value functions, this problem is solvable, and this will be the focus of the remainder of this section. The solution has two components: the first is an observation regarding the behavior of $E'_+$ and $E'_-$ under scaling. The second one is random variables, which have a stochastic representation in terms of a sum of random variables, where the number of terms follows a geometric distribution. We start with the first component. Suppose that $E'_+$ and $E'_-$ exists then we have  
\begin{align*}
\int_0^\infty E'_+(t)~dt=-\frac{2\mu_-}{\mu_++\mu_-}, 
\ \ \
\int_0^\infty E'_-(t)~dt=\frac{2\mu_+}{\mu_++\mu_-}.
\end{align*}
If we additionally assume that $E_+$ and $E_-$ are monotone functions, then the functions
\begin{align*}
    -\frac{\mu_++\mu_-}{2\mu_+}E'_+(t), 
    \ \ \ 
    \frac{\mu_++\mu_-}{2\mu_-}E'_-(t)
\end{align*}
are probability density functions.

The second component consists of the class of geometric divisible distributions. This generalizes the notion of geometric infinitely divisible, which was introduced by \cite{klebanov85}. The link between the geometric divisible distribution and non-stationary switch processes was noted in \cite{Bengtsson2}; this was, however, only for the case where $T_+\stackrel{d}{=}T_-$. We provide the following definition for completeness.  
\begin{definition}\label{rdiv}
Let $\nu_p$ be a geometric random variable with the probability mass function $p_{\nu_p} (k)=(1-p)^{k-1}p$ for $k=1,2...$ and  $\{V_k\}_{k\ge 1}$ a sequence of iid non-negative random variables independent of $\nu_p$. If a random variable has a stochastic representation 
\begin{align*}
W=\sum_{k=1}^{\nu_p} V_k,
\end{align*}
then $W$ follows a $r$-geometric divisible distribution with $r=1/p$ and said to belong to the class $GD(r)$, with the divisor $V$ and we write $F\in GD(r)$. 
\end{definition}
The notion that $V$ is a divisor of $W$ is natural since $V$ divides $W$ into smaller parts. From this definition, the density of $W$ can be expressed as $f_W=f_V^{\ast \nu_p}$. This form is unwieldy, and we therefore present the more useful expression in the Laplace domain 
\begin{align} \label{2}
   \Psi_W(s)=\frac{\frac{1}{r}\Psi_V(s)}{1-(1-\frac{1}{r})\Psi_V(s)}. 
\end{align} 
The main difference between the concepts of geometric divisibility and geometric infinite divisibility is for which $p\in (0,1)$ the stochastic representation of Definition \ref{rdiv} needs to hold. A random variable is said to be geometric infinitely divisible if this property is satisfied for \textit{all} $p\in(0,1)$, while for a random variable to be geometrically divisible, this property only needs to hold for \textit{one} specific $p\in(0,1)$. A simple consequence of this is that any geometric infinitely divisible distribution is also geometrically divisible. 

%
%

The next theorem connects non-stationary switch processes having monotone expected value functions with the class of geometrically divisible distributions. 

\begin{theorem} \label{th:1}
Let $D(t)$, $t\in\mathbb{R}$, be a non-stationary process with expected value functions such that $-E_+'(t)$ and $E_-'(t)$ are non-negative for $t>0$. Then, the 
time from a switch until it returns to that state has the the stochastic representation.
\begin{align*}
    T_- + T_+ \stackrel{d}{=} \sum_{k=1}^{\nu_\beta} X_k + \sum_{k=1}^{\nu_\alpha} Y_k 
\end{align*}
where $\nu_\alpha$ and $\nu_\beta$ are independent and geometrically distributed with the parameters
\begin{align*}
    \alpha=\frac{\mu_-}{\mu_++\mu_-}, \,\,
    \beta =\frac{\mu_+}{\mu_++\mu_-} 
\end{align*}
and the densities of $X_i$ and $Y_i$ are given by 
\begin{align*}
        f_X(t)=-\frac{1}{2\alpha}E_+'(t),\,\,
        f_Y(t)=\frac{1}{2\beta}E_-'(t),
\end{align*}
respectively. The sequences of $X_i$ and $Y_i$ are mutually independent sequences of iid random variables, which are also independent of $\nu_\alpha$ and $\nu_\beta$. 
\end{theorem}
\begin{proof}
    From Proposition~\ref{prop:invert}, the relation $\alpha=1-\beta$, and the following property of the Laplace transform $\mathcal{L}(h')(s)=s\mathcal{L}(h)-h(0)$ we have 
    \begin{align*}
        \Psi_+(s) \Psi_-(s)&=\frac{s\mathcal{L}(E_+)(s)-1}{s\mathcal{L}(E_-)(s)-1} \ \frac{s\mathcal{L}(E_-)(s)+1}{s\mathcal{L}(E_+)(s)+1},\\
        &=\frac{-\mathcal{L}(E_+')(s)}{2+\mathcal{L}
        (E_+')(s)} \ \frac{\mathcal{L}(E_-')(s)}{2-\mathcal{L}(E_-')(s)}, \\
        &=\frac{2\alpha \Psi_X(s) 2\beta \Psi_Y(s) }{(2-2\alpha\Psi_X(s))(2-2\beta\Psi_Y(s))}, \\
        &=\frac{\beta \Psi_X(s) }{1-(1-\beta)\Psi_X(s)} \ \frac{\alpha\Psi_Y(s)}{1-(1-\alpha)\Psi_Y(s)}. 
    \end{align*}
    The result follows since the equation above is the product of two equations, which are in the form of Equation (\ref{2}).
\end{proof}
While Theorem \ref{th:1} connects the distributions of the cycle length with functional properties of $E_+$ and $E_-$, it does not allow us to recover the distribution of $T_+$ and $T_-$. However, the next theorem allows us to recover the distribution of $T_+$ and $T_-$ from the expected value functions of the non-stationary switch process. 
\begin{theorem} \label{th:2}
Let $D(t)$, $t\in\mathbb{R}$ be a non-stationary switch process with the expected value functions $E_+$ and $E_-$. Then the following conditions are equivalent
\begin{itemize}
\item[\it i)] The functions $-E_+'(t)$ and $E_-'(t)$ are non-negative for $t>0$,
\vspace{2mm}
\item[\it ii)] 
$T_+$ and $T_-$ have the stochastic representations
    \begin{align*}
        T_+\stackrel{d}{=} X+\sum_{k=1}^{\nu_\alpha-1} Y_k, 
        \ \ \ \ 
        T_-\stackrel{d}{=} Y+\sum_{k=1}^{\nu_\beta-1} X_k, 
    \end{align*}
    \end{itemize}
    where $\nu_\alpha$ and $\nu_\beta$ are independent and geometrically distributed with the parameters
    \begin{align*}
        \alpha=\frac{\mu_-}{\mu_++\mu_-}&,\,\,\beta=\frac{\mu_+}{\mu_++\mu_-},
    \end{align*}
    and the densities of $X_i$ and $Y_i$ are given by 
    \begin{align*}
        f_X(t)=-\frac{1}{2\alpha}E_+'(t)&,\,\,
        f_Y(t)=\frac{1}{2\beta}E_-'(t). 
    \end{align*}
    respectively. The sequences of $X_i$ and $Y_i$ are mutually independent sequences of iid random variables, which are also independent of $\nu_\alpha$ and $\nu_\beta$
\end{theorem}

\begin{proof} $(i) \rightarrow (ii)$:
First, we note that if a random variable, $Z$,  has the stochastic representation $Z\stackrel{d}{=} X+\sum_{k=1}^{\nu_\alpha-1} Y_k$, then
\begin{align*}
    \Psi_Z(s)&=\Ex e^{s\left(X+\sum_{k=1}^{\nu_\alpha-1} Y_k\right)}
    =\Psi_X(s)\left(\Prob(\nu_\alpha=1)+\sum_{n=2}^\infty \Ex e^{s\sum_{k=1}^{n-1} Y_k} \Prob(\nu_\alpha=n)\right)\\
    &=\Psi_X(s)\sum_{n=1}^\infty \Psi_Y^{n-1}(s) (1-\alpha)^{n-1} \alpha
    =\frac{\alpha\Psi_X(s)}{1-\left(1-\alpha\right) \Psi_Y(s)}.
\end{align*}
    By Proposition~\ref{prop:invert}, Lemma \ref{lem:Ps}, the property $\mathcal{L}(h')(s)=s\mathcal{L}(h)-h(0)$, and $\alpha=1-\beta$ we have 
    \begin{align*}
        \Psi_+(s)&=\frac{s\mathcal{L}(E_+)(s)-1}{s\mathcal{L}(E_-)(s)-1}
        =\frac{-\mathcal{L}(E_+')(s)}{2-\mathcal{L}(E_-')(s)}
        =\frac{\alpha\Psi_X(s)}{1-\beta \Psi_Y(s)}
        =\frac{\alpha\Psi_X(s)}{1-(1-\alpha) \Psi_Y(s)}.
    \end{align*}
    The identical argument for $\Psi_-(s)$ is omitted. 
    \\
    $(ii)\rightarrow (i)$:
    From Lemma \ref{lem:Ps} and Corollary \ref{cor:Inv} we have that 
    \begin{align*}
        \mathcal{L}(E_+')(s)&=-2\Psi_+(s)\frac{1-\Psi_-(s)}{1-\Psi_+(s)\Psi_-(s)}
        \\ 
        -\frac{1}{2 \alpha }\mathcal{L}(E_+')(s)&=\frac{1}{\alpha}\frac{\alpha\Psi_X(s)}{1-\beta \Psi_Y(s)} 
        \frac{1-\frac{\beta\Psi_Y(s)}{1-\alpha \Psi_X(s)}}{1-\frac{\alpha\Psi_X(s)}{1-\beta \Psi_Y(s)} \frac{\beta\Psi_Y(s)}{1-\alpha \Psi_X(s)}}
        \\ 
        &= \Psi_X(s) \frac{1-\alpha \Psi_X(s) - \beta \Psi_Y(s)}{(1-\alpha \Psi_X(s))(1-\beta \Psi_Y(s))-\alpha \beta \Psi_X (s)\Psi_Y(s)}
        \\ 
        &=\Psi_X(s) \frac{1-\alpha \Psi_X(s) - \beta \Psi_Y(s)}{1-\alpha \Psi_X(s) - \beta \Psi_Y(s) +\alpha \beta \Psi_X(s) \Psi_Y(s) -\alpha \beta \Psi_X(s) \Psi_Y(s)} 
        \\
        &=\Psi_X(s).  
    \end{align*}
    Since $\Psi_X(s)$ is the Laplace transform of a density function, it follows that $(-2\alpha)^{-1}E'_+(t)$ will be non-negative and  $E_+(t)$ will therefore be a monotone function. The identical argument for $E_-(t)$ is omitted for brevity. 
\end{proof}

Theorem \ref{th:2} extend the results of \cite{Bengtsson2}, which showed that monotone expected values are equivalent to $T_+$ being $2-$geometric divisible when $T_+\stackrel{d}{=}T_-$. There is a straightforward interpretation for this case, with the geometric random variable being the number of switching attempts needed to switch, including the successful one. For the more general case when $T_+\stackrel{d}{\neq}T_-$ there is also a natural interpretation of Theorem \ref{th:2}. Specifically, one can view the time $X$ (or $Y$) as the random time until the first switching attempt is made. If this attempt is unsuccessful, $\nu_\alpha$ (or $\nu_\beta$) represents the number of failed switching attempts besides the first attempt, with the $Y$ (or $X$) being the random time between the attempts.

We conclude this section with a corollary on the implication of Theorem \ref{th:2} for the stationary switch process. 
\begin{corollary}
    Let $D(t)$, $t\geq0$, be a non-stationary switch process that satisfies the conditions of Theorem \ref{th:2}. Then the second derivative of the covariance function $R(t)$ of the stationary version $\tilde{D}(t)$, $t\in\mathbb{R}$ is of the form  
    \begin{align}
        R''(t)=\frac{4}{\mu_++\mu_-} \left(\alpha f_X(t) + (1-\alpha)f_Y\right(t)),
    \end{align}
    and is proportional to the probability density function of the mixture $\xi X+(1-\xi)Y$. Where $\xi$ is a Bernoulli variable with the parameter $\alpha$, independent of $X$ and $Y$. 
\end{corollary}
From this corollary, it is clear that even if the switching time distribution satisfies the condition of Theorem \ref{th:2}, they can not be recovered from the covariance function alone. This follows from the inability to decompose the density into a mixture of $X$ and $Y$ such that $X$ and $Y$ are unique.

\section{Examples} \label{sec:ex}
\subsection*{Asymmetry with a common divisor}
One of the simplest ways of inducing asymmetry is to let $T_+$ and $T_-$ be geometrically divisible with the same divisor $\tilde{T}$ but of different orders. The asymmetry is induced by modifying the single parameter, $\alpha\in(0,1)$ in the geometric summation, since $\beta=1-\alpha$. Hence, we have the following stochastic representation 
\begin{align*}
    T_+\stackrel{d}{=} \sum_{k=1}^{\nu_\alpha} \tilde{T}_k, \ \ \ 
        T_-\stackrel{d}{=} \sum_{k=1}^{\nu_\beta} \tilde{T}_k. 
\end{align*} 
Since this representation satisfies condition \textit{ii)} in Theorem \ref{th:2}, the expected value functions $E_+$ and $E_-$ will be monotone. Additionally, we have that $E_+/E_-=-\alpha/\beta$, meaning that they only differ by scaling. If the common divisor is exponential, this is similar to scaling $T_+$ and $T_-$ with some constant, as demonstrated in the next example.

\subsection*{Asymmetry with a common divisor trough scaling}
The asymmetry can be extended while maintaining the simplicity of having a common divisor, as in the previous example. The monotonicity property of the expected value functions $E_+$ and $E_-$ is preserved by considering a scaling of the form 
\begin{align*}
    T_+\stackrel{d}{=} b \tilde{T}_1+a\sum_{k=2}^{\nu_\alpha} \tilde{T}_k, \ \ \
        T_-\stackrel{d}{=} a \tilde{T}_1+b\sum_{k=2}^{\nu_\beta} \tilde{T}_k, \ \ a,b>0. 
\end{align*}
Since $\beta=1-\alpha$ this semi-parametric distribution has three numerical parameters $a>0$, $b>0$, $\alpha\in(0,1)$ and one functional parameter, which is the density $f$ of $\tilde{T}$. 
We observe that $\mu_+=\left(b+a(1/\alpha-1)\right) \mu$ and $\mu_-=\left(a+b(1/\beta-1)\right) \mu$, where $\mu=\Ex\tilde{T}$. Then the derivative of the expected value functions becomes 
\begin{align*}
E_+'(t)=-2\frac{a+b(1/\beta -1)}{a/\alpha+b/\beta}
f(t/b)/b, \ \ \ E_-'(t)=2\frac{b+a(1/\alpha -1)}{a/\alpha+b/\beta}f(t/a)/a.
\end{align*}

\subsection*{Non-monotonic expected value functions}
Consider a process with gamma-distributed switching times, with the scaling parameters $\theta_+=2$, $\theta_-=1$, and the shape parameters $k_+=2$ and $k_-=3$. The derivatives of the expected value functions can be evaluated via the Laplace transform using Corollary~\ref{cor:Inv}:
\begin{align*}
    \mathcal{L}(E'_+)(s)&=2 \,\frac{1-(1+s)^{3}}{(1+2s)^{2}(1+s)^3-1}=\frac{6+6s+2s^2}{7+19s+28s^2+16s^3+4s^4}, \\
   \mathcal{L}(E'_-)(s)&=2\, \frac{(1+2s)^{2}}{(1+2s)^{2}(1+s)^3-1}=\frac{8+8s}{7+19s+28s^2+16s^3+4s^4}. 
\end{align*}
The inverse Laplace transforms of these functions oscillate and provide an example where the conditions of Theorem \ref{th:2} are not met. A direct consequence is that the covariance function $R(t)$, of the corresponding stationary process will also oscillate due to the relation in  Proposition \ref{prop:rel}. 

\section{Conclusion}
\noindent
At the center of this paper are the characterizing functions of the two versions of the switch processes and their relationship to the switching time distributions. The main result shows that the monotonicity of the expected value functions is equivalent to the switching time distribution having a particular stochastic representation. This representation has a natural interpretation and allows us to recover the switching time distribution from the expected value functions. The recovery of the switching time distribution from the characterizing functions is one of the main motivations for studying these processes. We also show that the expected value functions should be the preferred characterization for this recovery. This is because using the covariance function of the stationary switch process leads to identifiability problems. These two observations have practical implications for the independent interval approximation framework, which approximates the excursions of Gaussian processes with a switch process. Applying these results to this framework is one of the two directions for future work. The other is to relax the independence assumptions between two switching time distributions. 

\section{Acknowledgement}
\noindent
Henrik Bengtsson acknowledges financial support from The Royal Physiographic Society in Lund. Krzysztof Podg\'orski and Henrik Bengtsson also acknowledge the financial support of the Swedish Research Council (VR) Grant DNR: 2020-05168.

\bibliographystyle{apalike}
\bibliography{References}
\newpage

\section*{Appendix. Auxiliary results and proofs}
\noindent
In this Appendix, we have collected some of the proofs for the paper along with some elementary and supporting results. Before these results are presented, it should be noted that we use $\ast$ for convolutions of functions and $\star$ for the convolution of probability distributions. Namely, $F\star F$ denotes the distribution function of the sum of two independent random variables that follow the distribution $F$. Lastly, the proofs are presented under the assumptions of Section \ref{sec:SP} unless stated otherwise. For convenience, the proofs and remarks of the appendix are grouped into three sections.

\subsection*{The Laplace transform and supporting results}
The Laplace transform is heavily used within renewal theory and in this paper. We therefore provide a short introduction of it here. Recall that for a function $h(t)$, $t\geq 0$, its Laplace transform $\mathcal L (h)(s)$, defined through 
\begin{align*}
\mathcal L (h)(s)=\int_0^\infty h(t)e^{-ts}~dt,
\end{align*}
for any $s$ such that the integral is finite. Furthermore, we have that for a positive half-line distribution given by a cdf $F$; its Laplace transform $\Psi$ is given by 
\begin{align*}
   \Psi(s)=\int_0^\infty e^{-ts} ~dF(t). 
\end{align*}
Additionally, there exists a well-known relation between the Laplace transform of the cumulative distribution function, $F$, and $\Psi$, namely
\begin{align} \label{cdf}
    \mathcal{L}(F)(s)=\Psi(s)/s.
\end{align}
which follows from the fact that $F$ is the convolution between the density $f$ and the Heaviside step function. Lastly, we want to highlight the following property of the Laplace transform of derivatives, which sees extensive use in this paper: 
\begin{align*}
    \mathcal{L}(h')(s)=s\mathcal{L}(h)(s)-h(0). 
\end{align*}
With this introduction to the Laplace transform and the assumptions in place, we are now ready to define the non-stationary switch process.

\begin{lemma}\label{fgsup}
Let $G$,$H$ be distribution functions on $[0,\infty)$ with corresponding  densities $g,h$. Then for each $t\ge 0$: 
\begin{align*}
    (g \ast h)(t) \le \min \left( \sup_{u>0}g(u) \cdot H(t), \ \sup_{u>0}h(u) \cdot G(t)\right).
\end{align*}
\end{lemma}
\begin{proof}
Since $g,h$ are densities on $[0,\infty)$, $g(t-x)=0$ for $x>t$, we have
\begin{align*}
    (g\ast h)(t)&=\int_{-\infty}^\infty g(t-x)h(x)dx=\int_{0}^t g(t-x)h(x)dx \\
    &\le  \int_{0}^t \sup_{u>0}g(u) h(x)dx=\sup_{u>0}g(u)H(t).
\end{align*}
We obtain the minimum since convolution is commutative,  
and we can freely choose to bound $g$ or $f$ with its supremum. 
\end{proof}

\begin{corollary}\label{corrbound}
Let $G$ be a distribution function on $[0,\infty)$, with density function $g$. Then for all $ n,m \in \mathbb{N}$ we have 
\begin{align*}
    g^{\ast (n+m)}(t) \le\sup_{u>0}g(u)^{\ast n} G(t)^{m\star}. 
\end{align*}
\end{corollary}
\begin{lemma} \label{eq:ConProd} 
Let $G$,$H$ be distribution functions on $[0,\infty)$, then we have 
\begin{align*}
    (H\star G)(t) \leq H(t) G(t)
\end{align*}
and for $n\in \mathbb{N}$, 
\begin{align*}
        H^{\star n}(t) \leq H^n(t).
\end{align*}
\end{lemma}
\begin{proof}
    Since $G$ and $H$ are distribution functions and thus monotone on $[0,\infty)$ we have  
    \begin{align*}
     (H\star G)(t)&=\int_{-\infty}^\infty H(t-u)~dG(u)=\int_0^t H(t-u)~dG(u) \le H(t)\int_0^t dG(u)=H(t)G(t). 
    \end{align*}
    We obtain the last equation from the associative property of convolutions and repeated application of the bound.  
\end{proof}

\subsection*{Proofs related to the non-stationary switch process}
After these supporting results, we account for the properties of the switch process $D(t)$ used across the paper. We start with the following result for the counting process $N(t)$ that counts the number of switches in $(0,t]$ for $t>0$ and those with the negative count, in $(t,0]$ for $t<0$.
Since we have 
\begin{align*}
D(t)=(-1)^{N(t)+(1-\delta)/2},
\end{align*}
it is clear that the distributional properties of $N(t)$ together with the initial state $\delta$ uniquely define the distribution of the process $D(t)$. The distribution of $N(t)$ on the positive half-line and conditionally on $\delta=\delta_0$, $\delta_0\in \{-1,1\}$ is the same as of $\{-N(-t)-1,t\ge 0\}$ conditionally on $\delta=-\delta_0$, hence it is sufficient to consider only the positive half-line.
\begin{lemma}
\label{cs}
Assume that $\mathcal T_0= (T_i^-,T_i^+)_{i \in \mathbb N}$ is a sequence of iid pairs of independent random variables, with the distribution function $F_{-}$ and $F_{+}$, respectively.  
The one-dimensional marginal distributions of $N(t)$   given $\delta$ and for $t\ge 0$ are given by
\begin{align} \label{rep}
\Prob{\left(N(t)=k \big | \delta\right)}&=\begin{cases}
\left(F_+ \star F_{-}\right)^{l\star}(t) -\left(F_+ \star F_{-}\right)^{l\star}\star F_\delta(t);& k=2l, \\
\left(F_+ \star F_{-}\right)^{l\star}\star F_{\delta}(t)-\left(F_+ \star F_{-}\right)^{l\star}(t);& k=2l+1,
\end{cases}
\end{align}
where $l\in \mathbb N\cup \{0\}$. 
\end{lemma}

\begin{proof}
We consider only the case of conditioning on $\delta=1$ as the opposite case can be obtained by the symmetry argument. 
We first note that for a positive $t$ and a non-negative integer $l$:
\begin{align*}
\Prob{\left(N(t)=0\big | \delta=1\right)}&=\Prob{(T^+_1>t)}=1- F_+(t),\\
\Prob{\left(N(t)=1\big | \delta=1\right)}&=\Prob{(T^+_1+T^-_1>t\ge T^+_1})\\
&=1-F_+\star F_-(t)-\left( 1- F_+(t)\right),\\
&= F_+(t)-F_+\star F_-(t),\\
\Prob{\left(N(t)=2l\big | \delta=1\right)}&=\Prob{\left(
\sum_{i=1}^{l+1}{T^+_i}+\sum_{i=1}^{l}{T^-_i}
>t
\ge 
\sum_{i=1}^{l}{T^+_i}+\sum_{i=1}^{l}{T^-_i}
\right)}\\
&=1-F_+^{(l+1)\star} \star F_-^{l\star }(t)-\left(1- F_+^{l\star}\star F_-^{l\star}(t)\right)\\
&= F_+^{l\star}*F_-^{l\star}(t)-F_+^{(l+1)\star}\star F_-^{l\star}(t),
\\
\Prob{\left(N(t)=2l+1\big | \delta=1\right)}&=\Prob{\left(
\sum_{i=1}^{l+1}{T^+_i}+\sum_{i=1}^{l+1}{T^-_i}
>t
\ge 
\sum_{i=1}^{l+1}{T^+_i}+\sum_{i=1}^{l}{T^-_i}
\right)}\\
&=F_+^{(l+1)\star}\star F_-^{l\star}(t)-F_+^{(l+1)\star} \star F_-^{(l+1)\star}(t).
\end{align*}
\end{proof}

With this lemma, we now present the proof of Proposition~\ref{OneDim}. 
\begin{proof}[Proof of Proposition~\ref{OneDim}]    Let us first consider the case of $t>0$ and $\delta =1$. 
Then by Equations (\ref{rep}), (\ref{cdf}) and Lemma~\ref{cs}:
\begin{align*}
P_+(t)&=\sum_{l=0}^\infty \Prob{(N(t)=2l)}=\left(\sum_{l=0}^\infty \left(F_+ \star F_-\right)^{l\star} \right) (t)-\left(\sum_{l=0}^\infty \left(F_+*F_-\right)^{l\star} \right)\star F_+ (t)\\
&=\mathcal L^{-1} \left(\frac{1-\Psi_+(s)}{s(1-\Psi_+(s)\Psi_-(s))}\right).
\end{align*}
For the case of $\delta=-1$, we have
\begin{align*}
P_-(t)&=\left(\sum_{l=0}^\infty \left(F_+ \star F_-\right)^{l\star} \right)\star  F_-(t)-\left(\sum_{l=1	}^\infty \left(F_+ \star F_-\right)^{l\star} \right)(t)\\
&=\mathcal L^{-1} \left(\frac{\Psi_-(s)\left(1-\Psi_+(s)\right)}{s(1-\Psi_+(s)\Psi_-(s))}\right).
\end{align*}
Moreover, for $\delta= 1$
\begin{align*}
E_{+}(t)&=2P_+(t)-1
\\
&=\mathcal L^{-1}\left(
 \frac{1-2\Psi_+(s)+\Psi_+(s)\Psi_-(s)}{s(1-\Psi_+(s)\Psi_-(s))} \right),
\end{align*}
while for $\delta=-1$ we obtain
\begin{align*}
E_{-}(t)&=
\mathcal L^{-1}\left(
 \frac{2\Psi_-(s)(1-\Psi_+(s))-1+\Psi_+(s)\Psi_-(s)}{s(1-\Psi_+(s)\Psi_-(s))}
\right)\\
&=\mathcal L^{-1}\left(
 \frac{-1+2\Psi_-(s)-\Psi_+(s)\Psi_-(s)}{s(1-\Psi_+(s)\Psi_-(s))}
\right).
\end{align*}
\end{proof}

The differentiability of $E_+$ and $E_-$ relies on the differentiability of $P_+$ and $P_-$ defined in Proposition~\ref{OneDim}. We provide the conditions for differentiability in the following lemma.
\begin{lemma}
\label{lem:Ps}
Assume that both distributions $F_+$ and $F_-$ are absolutely continuous with respect to the Lebesgue measure and have finite expectations, $\mu_+$ and $\mu_-$, respectively.
 The functions $P_+$ and $P_-$ satisfy the following conditions\vspace{2mm}
 \begin{itemize}
     \item[\it i)] $\displaystyle \lim_{t\rightarrow 0^+} P_{+}(t)=1$ and\, $\displaystyle \lim_{t\rightarrow 0^+} P_{-}(t)=0$\vspace{2mm},
     \item[\it ii)] $\displaystyle \lim_{t\rightarrow \infty} P_+(t)=\lim_{t\rightarrow \infty} P_-(t)=\frac{\mu_+}{\mu_++\mu_-}$ \vspace{2mm},
     \item[\it iii)] If  there exists $L\in \mathbb{N}$ for which $\sup_{t>0}(f_+*f_-)^{\ast L}(t)<\infty$, then $P_+$ and $P_-$ are differentiable and 
     \begin{align*}
P_+&=\left(\mathbb I_{[0,\infty)} - F_+ \right)\star \sum_{l=0}^\infty \left(F_+ \star F_-\right)^{l\star },\\
P_-&=F_- +\left(F_- -\mathbb I_{[0,\infty)}\right) \star\sum_{l=1}^\infty \left( F_+ \star F_-\right)^{l\star },\\
P_+'&= \sum_{l=0}^\infty \left(f_+^{*l  } * f_-^{*l } - f_+^{*(l+1)}* f_-^{l\star }\right),\\
P_-'&=f_+ +\sum_{l=1}^\infty \left(f_+^{*l  } * f_-^{*(l+1) } - f_+^{*l}* f_-^{l\star }\right).
\end{align*}
where the convergence of all series is locally uniform on $(0,\infty)$. 
 \end{itemize}
\end{lemma}
\begin{proof}
We have seen in the proof of Proposition~\ref{OneDim} that
\begin{align*}
P_+&=\mathbb I_{[0,\infty)} - F_+ +\sum_{l=1}^\infty \left(F_+^{l \star } \star F_-^{l\star } - F_+^{(l+1)\star}\star F_-^{l\star } \right)\\
&=\left(\mathbb I_{[0,\infty)}-F_+\right)\star U,
\end{align*}
where $U$ is the renewal measure defined by the renewal processes of the cycle length

\begin{align*}
U(A)=\delta_{\{0\}}+\sum_{l=1}^\infty \left(F_+ \star F_-\right)^{l\star }(A).
\end{align*}
from Lemma \ref{eq:ConProd} we have, 
\begin{align*}
    \lim_{t\rightarrow 0^+}\left(\mathbb I_{[0,\infty)} - F_+\right)(t)&=1,
    \\
    \lim_{t\rightarrow 0^+}\left|\sum_{l=1}^\infty \left(F_+^{l \star } \star F_-^{l\star } - F_+^{(l+1)\star}\star F_-^{l\star } \right)\right |(t)&\le \lim_{t\rightarrow 0^+}\frac{F_+\star F_-(t)}{1-F_+\star F_-(t)}(1+F_+(t)) =0,
\end{align*}
which establishes {\it i)} for $P_+$. The argument for $P_-$ is similar.

The line of the argument for {\it ii)} utilizes the Key Renewal Theorem, see p.86 in  \cite{theoryofpoint1}, which states that 
for an integrable function $g$ that vanishes on $(-\infty,0)$:
\begin{align*}
\lim_{t\rightarrow\infty}g\star U(t)=\frac{1}{\mu_++\mu_-}\int_0^\infty g(u)~du.
\end{align*}
For the limit at the infinity for $P_+$, one applies it to $g=\mathbb I_{[0,\infty)}-F_+$ since $P_+=g\star U$. The approach is analog for $P_-$

To prove the differentiability, we use the fact that $P_+$ and $P_-$ can be expressed as a series of differentiable terms; thus, from Theorem 7.17, \cite{RudinPoP}, it is sufficient to show local uniform convergence of the partial sums that would allow the exchange summation and derivation.
Thus for each $0< t_0 < t_1$, we need to show that, for each $\epsilon >0$, there exist an $N$ such that 
\begin{align*}
    \Big\vert \sum_{l=N+1}^\infty \left(f_+^{*l  } * f_-^{*l } - f_+^{*(l+1)}* f_-^{*l } \right)(t) \Big \vert  < \epsilon
\end{align*}
for all $t\in [t_0,t_1]$.

Assume first that the support of $f_+*f_-$ is the entire half-line $(0,\infty)$.
For given $t_0$ and $t_1$, let $t\in [t_0,t_1]$, and $N>L$. We note that $(F_+ \star F_-)(t_1)<1$.
Then for every $\epsilon>0$, by applying Corollary~\ref{corrbound} and Lemma \ref{eq:ConProd}, we have  
\begin{align*}
      \Big\vert 
     &\sum_{l=N+1}^\infty \left(f_+^{*l  } * f_-^{*l } - f_+^{*(l+1)}* f_-^{l\star } \right)(t)
      \Big\vert  
      \\ 
      &\le \left(f_+*f_-\right)^{*L}*   \sum_{n=1}^\infty  \left( \left(f_+*f_-\right)^{\ast (N-L+n)} +\left(f_+*f_-\right)^{\ast (N-L+n)}*f_+\right)(t) 
      \\
     &\le \sup_{u>0} \left(f_+*f_-\right)^{\ast L}(u) \cdot \sum_{n=1}^\infty  \left( \left(F_+ \star F_-\right)^{ (N-L+n)\star} +\left(F_+\star F_-\right)^{ (N-L+n)\star}\star F_+\right)(t)
     \\
     &\le \sup_{u>0} \left(f_+*f_-\right)^{\ast L}(u) \cdot \sum_{n=1}^\infty   \left(F_+ \star F_-\right)^{ (N-L+n)}(t_1) \left( 1+ F_+(t_1)\right)
     \\
     &=\sup_{u>0} \left(f_+*f_-\right)^{\ast L}(u) \cdot  \frac{\left(F_+ \star F_-(t_1)\right)^{ (N-L+1)} }{1-F_+ \star F_-(t_1)}\left( 1+ F_+(t_1)\right),
\end{align*}
where the last term can be made smaller than $\epsilon$ by taking $N$ sufficiently large so that $\left(F_+ \star F_-(t_1)\right)^{ (N-L+1)}$ is small enough. 
This proves uniform convergence of $P'_+$ and $P'_-$ on compact subsets of $\mathbb{R}^+$.
The proof can be extended for when the support of $f_+*f_-$ is not the entire positive real line; this is technical and, therefore, omitted in this paper. However, the approach is to note that convolution extends the support. In essence, there is some $k$ chosen large enough such that our $t_1$ falls within the support and $F_+\star F_-^{k \star}(t_1)<1$. Then, this is used as the size of the groupings of the subsequent terms in the series. 
\end{proof}

\begin{proof}[Proof of Proposition \ref{prop:invert}]
Form the observation that $\mathcal{L}(E_\pm)(s)=2\mathcal{L}(P_\pm)(s)-1/s$ we have the following system of equations
\begin{align*}
          \begin{cases}
          \mathcal{L}(E_+)(s)+\frac{1}{s}&=2\mathcal{L}(P_+)(s) \\
          \mathcal{L}(E_-)(s)+\frac{1}{s}&=2\mathcal{L}(P_-)(s).
          \end{cases}
\end{align*}
By using the fact from Proposition \ref{OneDim} that $\mathcal{L}(P_-)(s)/\mathcal{L}(P_+)(s)=\Psi_-(s)$ and the first equation of the same proposition we obtain 
\begin{align*}
          \begin{cases}
          \mathcal{L}(E_+)(s)+\frac{1}{s}&=\frac{2}{s} \frac{1-\Psi_+(s)}{1-\Psi_+(s)\Psi_-(s)} \\
          \frac{s\mathcal{L}(E_-)(s)+1}{s\mathcal{L}(E_+)(s) +1}&=\Psi_-(s).
          \end{cases}
\end{align*}
Now, we can solve the top equation by substitution, and we have 
\begin{align*}
    (s\mathcal{L}(E_+)(s)+1) \left( 1 - \Psi_+(s) \frac{s\mathcal{L}(E_-)(s)+1}{s\mathcal{L}(E_+)(s) +1}  \right) &=2(1-\Psi_+(s))
    \\
    (s\mathcal{L}(E_+)(s)-1) - \Psi_+(s)\left( s\mathcal{L}(E_-)(s) -1 \right)&=0, 
\end{align*}
which leads to 
\begin{align*}
    \begin{cases}
           \frac{s\mathcal{L}(E_+)(s)-1}{s\mathcal{L}(E_-)(s)-1}&=\Psi_+(s) \\
           \frac{s\mathcal{L}(E_-)(s)+1}{s\mathcal{L}(E_+)(s)+1}&=\Psi_-(s). 
    \end{cases} 
\end{align*}
\end{proof}
\subsection*{Proofs related to the stationary switch process}
Before moving to the results on the stationary switch process, we need to define some well-known objects in renewal theory in the context of this paper. Define $S_L(t)$ and $S_R(t)$ as the instants of the most recent switch prior to $t$  and the first switch after $t$, respectively. 
The residual time since the last switch is denoted by $B(t)=t-S_L(t)$, and the excess time is $A(t)=S_R(t)-t$. 
\begin{lemma}
\label{compl}
The joint distribution of $A(t),B(t),N(t)$ is given through the following conditional probabilities for $u\in [0,t]$ and $v>0$:
\begin{align*}
 \Prob{\left(A(t)> v\big |B(t)=u, N(t)=k, \delta \right)}&
 =\begin{cases}
\bar F_\delta(v+u)/\bar F_\delta(u);& k=2l ,\\
\bar F_{-\delta}(v+u)/\bar F_{-\delta}(u);& k=2l+1,
\end{cases}
\end{align*}
\begin{align*}
 \Prob{\left(S_L(t)\le u\big | N(t)=k, \delta \right)}&
 =\begin{cases}
F_\delta^{l\star}\star F_{-\delta}^{l\star}(u)\cdot \overline{\left(F_\delta^{l\star }\star F_{-\delta}^{l\star }\right)_u\star F_\delta}(t) ;& k=2l ;\\
F_\delta^{(l+1)\star }\star F_{-\delta}^{l\star }(u)\cdot \overline{\left(F_{\delta}^{(l+1)\star} \star F_{-\delta}^{l\star}\right)_u\star F_{-\delta}}(t) ;& k=2l+1;
\end{cases}
 \end{align*}
A probability measure represented by $F$ which is restricted to $[0,u]$ is denoted by $F_u$ i.e., $F_u(x)=F(x)/F(u)$ for $x\in [0,u]$ and one for $x>u$. 

Further, the process $\Delta N(s)=N(s+S_R(t))-N(S_R(t))$, $s>0$, depends on the triplet  $\left(S_L(t),S_R(t),N(t)\right)$ only through $\delta_0=(-1)^{N(t)+1}$. Hence, for a given $\delta=\delta_0$, it has the same distribution (as the entire process) as the process $N(s)$, $s>0$.
\end{lemma}
\begin{proof}
For $B(t)=u,N(t)=k,\delta=1$, we have $\sum_{i=1}^{l}{T^+_i}+\sum_{i=1}^{l}{T^-_i}=t-u$, and  $t-u + T^+_{l+1}>t$, i.e. $T^+_{l+1}>u$. 
Given these conditions, $A(t)>v$ is equivalent to $t-u+T^+_{l+1}>t+v$, i.e. $T^+_{l+1}>u+v$.
Thus 
\begin{align*}
 \Prob{\left(A(t)> v\big |B(t)=u, N(t)=k, \delta=1 \right)}&
 =\begin{cases}
 \Prob(T^+_{l+1}-u>v\big | T^+_{l+1}>u)&; k=2l,\\
 \Prob(T^-_{l+1}-u>v\big | T^-_{l+1}>u)&; k=2l+1,
 \end{cases}\\
 &=\begin{cases}
 {\bar F_+(u+v)}/{\bar F_+(u)}&; k=2l,\\
 {\bar F_-(u+v)}/{\bar F_-(u)}&; k=2l+1,
 \end{cases}
 \end{align*}
The case of $\delta=-1$ is symmetric. 

Further,
 \begin{align*}
 \Prob{\left( B(t)> u\big |  N(t)=2l, \delta=1 \right)}&=\Prob{\left(
t-u > \sum_{i=1}^{l}{T^+_i}+\sum_{i=1}^{l}{T^-_i} > t- T_{i}^+
\right)},\\
&=\int_0^{t-u} \bar F_+(t-p) d F_+^{l\star}\star F_-^{l\star }(p)\\
&=F_+^{l\star} \star F_-^{l\star }(u)\cdot \overline{\left(F_+^{l\star }\star F_-^{l\star }\right)_u\star F_+}(t),\\
\Prob{\left( B(t)\le u\big |  N(t)=2l+1, \delta=1 \right)}&=\Prob{\left(
\sum_{i=1}^{l+1}{T^+_i}+\sum_{i=1}^{l}{T^-_i}
\le u
, 
\sum_{i=1}^{l+1}{T^+_i}+\sum_{i=1}^{l+1}{T^-_i}>t
\right)},\\
&=\int_0^u \bar F_-(t-p) d F_+^{(l+1)\star }\star F_-^{l\star }(p)\\
&=F_+^{(l+1)\star }\star F_-^{l\star }(u)\cdot \overline{\left(F_+^{(l+1)\star }\star F_-^{l\star }\right)_u\star F_-}(t).
 \end{align*}
The second part of the result follows from the mutual independence of $T_i^+$'s and $T_i^-$'s. 
\end{proof}
The following proof uses the key renewal theorem or, more precisely, one of its consequences, the limiting behavior of alternating processes. 

\begin{proof}[The proof of Proposition~\ref{stat}]
Let $D(t)$ be the regular switch process starting at time zero. 
Since the distribution of $N(s)$ is the same as $\Delta N(s)$ in Lemma~\ref{compl}, the stationary distribution can be obtained by finding the distributional limit of the triple $(A(t),B(t),D(t))$. The limit of the triple can be obtained from the Key Renewal Theorem.

Consider the renewal process defined via $X_n=T_n^++T_n^-$ and let $S_n$ be the cumulative sum of $X_n$. 
By defining $Z_n=T_n^+$, we obtain an alternating process that is ON over $[S_n, S_n+T_{n+1}^+)$, $n\in \mathbb N$, and OFF, otherwise. 
Thus, the limiting distribution of probability that at time $t$ the process is in the ON state is $\Ex Z_1/\Ex X_1$, see Theorem~3.4.4 in \cite{Ross}.  
This essentially proves the limiting distribution of $D(t)$. 

To show the joint limiting distribution of $A(t)$ and $B(t)$ given $D(t)=1$, we note that for an alternating process defined by $Z_n=T_n^+\wedge( T_n^+-b) \wedge a$   the system is ON if and only if $D(t)=1, A(t)\le a, B(t)>b$. 
We note that $\Ex Z_n=\int_0^a \bar F_+(u+b)~du$ and thus the limiting probability for $\Prob{(\delta(t)=1, A(t)\le a, B(t)>b)}$ is equal to $\int_0^a \bar F_+(u+b)~du/(\mu_++\mu_-)$, which demonstrates the distributional form of $(\delta, A, B)$, see Proposition~3.4.5, \cite{Ross}. 

The expected value of the stationary switch process $\tilde D$ is given by $\Ex \delta=(\mu_+-\mu_-)/(\mu_++\mu_-)$.\
To find $\Prob{(\tilde D(t)=1)|\tilde D(0)=1})$, we follow the same path as in the proof of Proposition~\ref{OneDim}. 
If $\tilde N(t)$ denotes the number of switches between $0$ and $t$ for the stationary switch process, then its distribution at time $t$ is given by
\begin{multline*}
\Prob{\left(\tilde N(t)=k \big | \tilde D(0)= \delta\right)}=\\
=
\begin{cases}
1-F_{A,\delta}(t); k=0 \\
 \left(F_+ \star  F_{-}\right)^{(k/2-1)\star }\star F_{-\delta}\star  F_{A,\delta}(t)-\left(F_+ \star  F_{-}\right)^{(k/2)\star}\star F_{A,\delta}(t);& k>0 \mbox{ is even},\\
\left(F_+ \star  F_{-}\right)^{(k-1)/2\star}\star F_{A,\delta}(t)-\left(F_+ \star F_{-}\right)^{(k-1)/2\star}\star F_{-\delta}\star F_{A,\delta}(t);& k \mbox{ is odd},
\end{cases}
\end{multline*}
where $F_{A,\delta}$  is the cdf of the initial delay given in Remark~\ref{AB}.
The proof is a slight modification of the one in Theorem~\ref{cs} in which one needs to account for the delay $A$ distribution that is different from $T_1^+$, see Remark~\ref{AB}.

Since 
\begin{align} \label{repd}
\tilde D(t)=(-1)^{N_s(t)+(1-\delta)/2}
\end{align} 
we can get, in the same way as in the proof of Proposition~\ref{OneDim}, that
\begin{align*}
\tilde P_{+}(t)&= \Prob{(\tilde D(t)=1)|\tilde D(0)=1})
=\sum_{l=0}^\infty \Prob(\tilde N(t)=2l \big | \tilde D(0)=1)
\\
 &=1-F_{A,+}(t) +  \sum_{l=1}^\infty \left(F_+ \star F_{-}\right)^{(l-1)\star }\star F_{-}\star F_{A,+}(t)
  - \sum_{l=1}^\infty \left(F_+ \star F_{-}\right)^{l\star }\star F_{A,+}(t)\\
 &=\mathcal L^{-1} \left(\frac 1s \left(1-\frac{1-\Psi_{A,+}(s)}{(1-\Psi_+(s)\Psi_-(s))}\right)\right)(t)
\end{align*}
 and 
\begin{align*}
    \tilde P_{-}(t)=\Prob{(\tilde D(t)=1)|\tilde D(0)=-1})=\mathcal L^{-1} \left(\frac{\Psi_{A,-}(s)\left(1-\Psi_+(s)\right)}{s(1-\Psi_+(s)\Psi_-(s))}\right).
 \end{align*}
Then the final formula for the Laplace transform of $P_\delta(t)$ follows from the form of $\Psi_{A,\delta}(s)$ that is given in Remark~\ref{AB}. 

Finally, for the covariance function, we note that
\begin{align*}
R(t)&=\Ex(\tilde D(0)\tilde D(t))-\frac{(\mu_+-\mu_-)^2}{(\mu_++\mu_-)^2}\\
&= \left(2\tilde P_{+}(t)-1\right)\frac{\mu_+}{\mu_++\mu_-}-\left(2\tilde P_{-}(t)-1\right)\frac{\mu_-}{\mu_++\mu_-}-\frac{(\mu_+-\mu_-)^2}{(\mu_++\mu_-)^2}\\
 &=\frac{2}{\mu_++\mu_-}
 \left(
 \tilde P_{+}(t)\mu_+-\tilde P_{-}(t)\mu_-+
\mu_+ 
 \frac{\mu_--\mu_+}{\mu_++\mu_-}
 \right)\\
\end{align*}
The Laplace transform of $R$ can now be easily obtained by combining all the above results. 
\end{proof}

\begin{proof}[Proof of Proposition \ref{prop:rel}]
    The first equations we obtain from Proposition \ref{stat} and \ref{OneDim}
    \begin{align*}
        \mathcal{L}(\tilde{P_+'})(s)&=s\mathcal{L}(\tilde{P_+})(s)-\Prob (X_+(0)=1 \vert X_+(0)=1),
        \\ %
        &=1-\frac 1 {\mu_+ s} \frac{\left(1-\Psi_+(s)\right)\left(1-\Psi_-(s)\right)}{1-\Psi_+(s)\Psi_-(s)} -1,
        \\ %
        &= \frac 1 {\mu_+ s} \frac{\left(1-\Psi_+(s)\right)\left(1-\Psi_-(s)\right)}{1-\Psi_+(s)\Psi_-(s)},
        \\%
        &=\frac{1}{\mu_+ }\left( \frac{1-\Psi_+(s)}{s(1-\Psi_+(s)\Psi_-(s))}-\frac{\Psi_-(s)(1-\Psi_+(s))}{s(1-\Psi_+(s)\Psi_-(s))}\right), 
        \\%
        &=-\left( \frac{\mathcal{L}(P_+)(s)-\mathcal{L}(P_-)(s)}{\mu_+}\right).
    \end{align*}
    The proof is very similar for $\mathcal{L}(\tilde{P_+'})(s)$ and observing that from the construction of $\tilde{P_+'}$ if follows that $\Prob (X_-(0)=1 \vert X_-(0)=-1)=0$. 
    The second equation follows from substituting in the first and that $E_+(t)=2\Prob (X_+(t)=1)-1$.
\end{proof}



\end{document}